\newtheorem{defn}{Definition}

\newtheorem{thm}{Theorem}
\newtheorem{prop}{Proposition}
\newtheorem{lem}{Lemma}
\newtheorem{rem}{Remark}

\newcommand{\ung}{{\mbox{$\underline{g}$}}}

\newcommand{\uqsl}{{\mbox{$U_{q}(sl_{n+1})$}}}
\newcommand{\uqgl}{{\mbox{$U_{q}(gl_{n+1})$}}}

\newcommand{\uqsla}{{\mbox{$U_{q}({\hat {sl}}_{n+1})$}}}
\newcommand{\uqslp}{{\mbox{$L_{q}({ {sl}}_{n+1})$}}}

\newcommand{\unv}{{\mbox{$\bf{v}$}}}
\newcommand{\una}{{\mbox{$\bf{a}$}}}

\newcommand{\uns}{{\mbox{$\bf{s}$}}}

\newcommand{\unj}{{\mbox{$\bf{j}$}}}

\newcommand{\bed}{\begin{defn}}
\newcommand{\bet}{\begin{thm}}
\newcommand{\bep}{\begin{prop}}
\newcommand{\bel}{\begin{lem}}
\newcommand{\brk}{\begin{rem}}

\newcommand{\ot}{{\mbox{$\otimes$}}}

\newcommand{\eed}{\end{defn}}
\newcommand{\eet}{\end{thm}}
\newcommand{\eep}{\end{prop}}
\newcommand{\eel}{\end{lem}}
\newcommand{\erk}{\end{rem}}

\newtheorem{cor}{Corollary}
\newcommand{\becor}{\begin{cor}}
\newcommand{\eecor}{\end{cor}}

\newcommand{\ci}{{\mbox{$C\!\!\!\!I$}}}
\newcommand{\sg}{{\mbox{$S_{\ell}$}}}
\newcommand{\hm}{{\mbox{$H_{\ell}(q^2)$}}}
\newcommand{\hma}{{\mbox{${\hat {H}}_{\ell}(q^2)$}}}

\newcommand{\hmi}{{\mbox{$H_{\ell_i}(q^2)$}}}
\newcommand{\hmai}{{\mbox{${\hat {H}}_{\ell_i}(q^2)$}}}

\newcommand{\hmo}{{\mbox{$H_{\ell_1}(q^2)$}}}

\newcommand{\hmt}{{\mbox{$H_{\ell_2}(q^2)$}}}

\newcommand{\vom}{{\mbox{$V_0^{\ot{\ell}}$}}}

\newcommand{\vm}{{\mbox{$V^{\ot{\ell}}$}}}

\newcommand{\uqg}{{\mbox{$U_q(\ung)$}}}

\newcommand{\cit}{{\mbox{$\ci^{\times}$}}}

\documentstyle[11pt]{article}
\batchmode
\begin{document}
\centerline{ {\Large Quantum Affine Algebras and Affine Hecke Algebras}}
\vskip 24pt
\centerline{\large Vyjayanthi Chari{\footnote{ Both authors were partially
supported by the NSF, DMS-9207701.}} and Andrew Pressley}
\vskip 36pt
\centerline{ August 30, 1993}
\vskip 36pt
\section{Introduction}

One of the most beautiful results from the classical period of the
representation theory of Lie groups is the correspondence, due to Frobenius and
Schur, between the representations of symmetric groups and those of general or
special linear groups. If $V_0$\/ is the natural irreducible
$(n+1)$--dimensional representation of $SL_{n+1}(\ci)$, the symmetric group
$S_\ell$\/ acts on $V_0^{\otimes\ell}$\/ by permuting the factors. This action
obviously commutes with the action of $SL_{n+1}(\ci)$. It follows that one may
associate to any right \sg--module $M$\/ a representation of $SL_{n+1}(\ci)$\/,
namely
$${\cal FS}(M) = M\ot_{\sg}V_0^{\otimes\ell},$$
the action of $SL_{n+1}(\ci)$\/ on ${\cal FS}(M)$\/ being induced by its
natural action on $V_0^{\otimes\ell}$. The main result of the Frobenius--Schur
theory is that, if $\ell\le n$, the assignment $M\to {\cal FS}(M)$\/ defines an
equivalence from the category of finite--dimensional representations of \sg\/
to the category of finite--dimensional representations of $SL_{n+1}(\ci)$, all
of whose irreducible components occur in $V_0^{\otimes\ell}$.

Around 1985, Drinfeld and Jimbo independently introduced a family of Hopf
algebras \uqg , depending on a parameter $q\in\ci^{\times}$,  associated to any
symmetrizable Kac--Moody algebra \ung .
Assuming that $q$\/ is not a root of unity, Jimbo [7] proved an analogue of the
Frobenius--Schur correspondence in which $SL_{n+1}(\ci)$\/ is replaced by \uqsl
, $V_0$\/ by the natural $(n+1)$--dimensional irreducible representation $V$\/
of \uqsl, and \sg\/ by its Hecke algebra \hm .

In [5], Drinfeld announced an analogue of the Frobenius--Schur theory for the
Yangian $Y(sl_{n+1})$, which is a  \lq\lq deformation" of the universal
enveloping algebra of the Lie algebra of polynomial maps $\ci\to sl_{n+1}$. The
role of \sg\/ in this theory is played by the degenerate affine Hecke algebra
$\Lambda_\ell$, an algebra whose defining relations are obtained from those of
the affine Hecke algebra \hma\/ by letting $q\to 1$\/ in a certain non--trivial
fashion.

In the same paper, Drinfeld conjectured that there should be an analogue of the
Frobenius--Schur theory  relating the quantum affine algebra \uqsla\/ and \hma
{}.
In this paper, we construct a functor from the category of finite--dimensional
\hma--modules to the category of finite--dimensional \uqsla--modules $W$\/ of
\lq type 1' (a mild spectral condition) with the property that every
irreducible \uqsl--type which occurs in $W$\/ also occurs in
$V^{\otimes\ell}$\/ (we assume that $q$\/ is not a root of unity). We prove
that this functor is an equivalence if $\ell\le n$. Drinfeld's theory can be
obtained from ours by taking a suitable limit $q\to 1$. Related results were
obtained by Cherednik in [4].

We give a precise description of our functor at the level of irreducible
representations, using the known parametrizations of such representations of
\uqsla\/ and of \hma. Namely, in [2], [3] we showed that the
finite--dimensional irreducible \uqsla--modules of type 1  are in one to one
correspondence
with $n$--tuples of monic polynomials in one variable.  On the other hand,
Zelevinsky [13] and Rogawski [12] have given a one to one correspondence
between the finite--dimensional irreducible \hma--modules and the set of
(unordered) collections of \lq segments\rq\/ of complex numbers, the sum of
whose lengths is $\ell$. (A segment of length $k$\/ is a $k$--tuple of the form
$(a, q^2a,\ldots ,q^{2k-2}a)$, for some $a\in\cit$.) We compute explicitly the
$n$--tuple of polynomials associated under our functor to any such collection
of segments.

The affine Lie algebra $\hat{sl}_{n+1}$\/ is a central extension, with
one-dimensional centre, of the Lie algebra of Laurent polynomial maps
$f:\cit\to sl_{n+1}$. An obvious way to construct representations of
$\hat{sl}_{n+1}$\/ is to pull back a representation of $sl_{n+1}$\/ by the
one-parameter family of homomorphisms $ev_a^0:\hat{sl}_{n+1}\to sl_{n+1}$\/
which annihilate the centre and evaluate the maps $f$\/ at $a\in\ci^\times$. In
[7], Jimbo defined a one-parameter family of algebra homomorphisms
$ev_a:\uqsla\to\uqsl$\/ which are quantum analogues of the $ev_a^0$\/
(actually, $ev_a$\/ takes values in an `enlargement' of \uqsl). On the other
hand, in [4] Cherednik defined a one-parameter family of homomorphisms
$\tilde{ev}_a:\hma\to\hm$\/ which are the identity on $\hm\subset\hma$. Pulling
back representations of \uqsl\/ (resp. \hm) under $ev_a$\/ (resp.
$\tilde{ev}_a$) gives a one-parameter family of representations of \uqsla\/
(resp. \hma). We show that these `evaluation' representations correspond to

 each other under our functor.
\vskip 12pt\noindent{\bf Acknowledgement} We would like to thank I. V.
Cherednik for several discussions related to this work.

\section{Quantum Kac--Moody algebras}
Let $A = (a_{ij})$\/ be a symmetric generalized Cartan matrix, where the
indices $i$, $j$\/ lie in some finite set $I$. Thus, $a_{ij}\in{\bf Z}$,
$a_{ii} =2$, and $a_{ij}\le 0$\/ if $i\ne j$. To $A$\/ one can associate a
Kac--Moody Lie algebra $\ung(A)$\/ (see [8]).

Let $q$\/ be a non--zero complex number, assumed throughout this paper not to
be a root of unity. For $n, r\in {\bf N}$, $n\ge r$, define
\begin{eqnarray*} [n]_{q} &=&\frac{q^n -q^{-n}}{q -q^{-1}},\\
{}\left[{n\atop r}\right]_{q}&=&\frac{[n]_{q}[n-1]_{q}\ldots {}[n-r+1]_{q}}
{[r]_{q}[r-1]_{q}\ldots [1]_{q}}.\end{eqnarray*}

\subsection{}\bed The quantum Kac--Moody algebra $U_q(\ung(A))$\/ associated to
a symmetric generalized Cartan matrix $A= (a_{ij})_{i,j\in I}$\/ is the unital
associative algebra over \ci\/ with generators $x_i^{{}\pm{}}$, $k_i^{{}\pm
1}$\/ ($i\in I$) and the following defining relations:
\begin{eqnarray*}k_ik_i^{- 1} =&1& =  k_i^{ -1}k_i,\\
k_ik_j&=&k_jk_i\;,\\
k_ix_j^{{}\pm{}}k_i^{-1}&=& q^{{}\pm a_{ij}}x_j^{{}\pm{}}\; ,\\
{}[x_i^+ , x_j^-]&=&\delta_{ij}\frac{k_i -k_i^{-1}}{q
-q^{-1}}\;,\end{eqnarray*}
$$\sum_{r=0}^{1-a_{ij}}
 (-1)^r\left[{1-a_{ij}\atop r}\right]_{q}\;
(x_i^{{}\pm{}})^rx_j^{{}\pm{}}(x_i^{{}\pm{}})^{1-a_{ij}-r}= 0\;,\;i\ne j .$$
\eed

It is well--known that $U_q(\ung(A))$\/ is a Hopf algebra with comultiplication
$\Delta$\/ given on generators by
\begin{eqnarray*}
\Delta(k_i^{{}\pm 1}) &=& k_i^{{}\pm 1}\ot k_i^{{}\pm 1},\\
\Delta(x_i^+) &= &x_i^+\ot k_i +1\ot x_i^+,\\
\Delta(x_i^-) &=& x_i^-\ot 1 + k_i^{-1}\ot x_i^- \end{eqnarray*}
(we shall not need the formulas for the counit  and antipode of
$U_q(\ung(A))$).

\subsection{} By a representation of a quantum Kac--Moody algebra
$U_q(\ung(A))$\/ we shall mean a left
$U_q(\ung(A))$--module. A representation $W$\/ is  said to be of type 1 if
$$W=\bigoplus_{{\bf\mu}\in {\bf Z^I}}W_{{\bf\mu}}, $$
where $W_{\mu} =\{w\in W|k_i.w = q^{\mu(i)} w\}$. If $W_{\mu}$\/ is non--zero,
then $W_{\mu}$\/ is called the weight space of $W$\/ with weight $\mu$.
Restricting consideration to type 1 representations results in no essential
loss of generality, for any finite--dimensional irreducible representation can
be obtained by twisting a type 1 representation with a suitable automorphism of
$U_q(\ung(A))$\/ (cf.  [10]).

\subsection{} Assume that ${\rm dim}(\ung(A))  < \infty$. A  representation
$W$\/ of $U_q(\ung(A))$\/ is said to be highest weight with highest weight
$\lambda\in {\bf Z^I}$\/ if $W$\/ is generated as a $U_q(\ung(A))$--module by
an element $w_{\lambda}$\/  satisfying
$$x_i^+. w_{\lambda} =0,\;\;k_i. w_{\lambda} =q^{\lambda(i)}w_{\lambda},$$
for all $i\in I$.

A weight $\lambda\in{\bf Z^I}$\/ is said to be dominant  if $\lambda(i)$\/ is
non--negative  for all $i\in I$.

\bep{{\rm ([10])}} Assume that ${\rm dim}\;(\ung(A)) < \infty$.

\noindent (i) Every finite--dimensional $U_q(\ung(A))$--module is completely
reducible.

\noindent (ii) Every irreducible finite--dimensional $U_q(\ung(A))$--module of
type 1 is  highest weight. Assigning to such a representation its highest
weight defines a one to one correspondence between the set of isomorphism
classes of finite--dimensional irreducible representations of type 1 and the
set of dominant weights.

\noindent (iii) The finite--dimensional irreducible $U_q(\ung(A))$--module
$V(\lambda)$\/ of type 1 and highest weight $\lambda$\/ has the same character
(in particular, the same dimension) as the irreducible $\ung(A)$--module of the
same highest weight.

\noindent (iv) The multiplicities of the irreducible components in a tensor
product $V(\lambda)\ot V(\mu)$\/ of irreducible finite--dimensional
$U_q(\ung(A))$--modules is the same as in the tensor product of the irreducible
$\ung(A)$--modules of the same highest weights.
$\Box$\eep

\subsection{} The case of most interest to us is when $A$\/ is the matrix
$$\left(\begin{array}{ccccccc}2&-1&0&0&\cdots&0&-1\\
-1&2&-1&0&\cdots&0&0\\
0&-1&2&-1&\cdots&0&0\\
\vdots&\vdots&\vdots&\ddots&\vdots&\vdots&\vdots\\
0&0&0&\cdots&-1&2&-1\\
-1&0&0&\cdots&0&-1&2
\end{array}\right),$$
where $i,j\in\{0,1,\ldots , n\}$. Then $\ung(A)$\/ is the affine Lie algebra
${\hat{sl}}_{n+1}$. Fix a square root $q^{1/2}$\/ of $q$. For any elements $a$,
$b$\/ of an associative algebra over $\ci$, set
$$[a,b]_{q^{1/2}}=q^{1/2}ab-q^{-1/2}ba.$$
Since $a_{ij}=0$\/ or $-1$\/ if $i\ne j$, the quantized Serre relations in
\uqsla\/ can be written
\begin{eqnarray*}
{}[x_i^{{}\pm{}},x_j^{{}\pm{}}]&=&0\;\;{\rm if}\;i-j\ne 0,\;{{}\pm 1}\;({\rm
mod }\;n),\\
{}[x_i^{{}\pm{}},[x_j^{{}\pm{}},x_i^{{}\pm{}}]_{q^{1/2}}]_{q^{1/2}}&=&0\;\;
{\rm if}\;i-j={{}\pm 1}\;({\rm mod} \;n).
\end{eqnarray*}

Deleting the $0^{th}$\/ row and column of $A$\/ gives the Cartan matrix of
$sl_{n+1}$. Thus, there is a natural Hopf algebra  homomorphism from \uqsl\/ to
\uqsla ; this homomorphism is injective (this follows from Proposition 5.4
below).

If $\ung(A) = sl_{n+1}$, then $I =\{1,\ldots ,n\}$\/ and so weights are
identified with $n$--tuples of integers.  It is useful to introduce the weights
$\epsilon_i$, for $1\le i\le n$, defined by
$$\epsilon_i(j) =\left\{ \begin{array}{rl}
 -1& {\rm if}\; j= i-1,\\
1& {\rm if}\; j =i,\\
 0&{\rm otherwise}.\end{array}\right.$$
Note that $\sum_{i=1}^{n+1}\epsilon_i = 0$.

Set $\alpha_i=\epsilon_i-\epsilon_{i+1}$. If $\lambda$, $\mu\in{\bf Z^I}$, we
write $\lambda\ge\mu$ if $\lambda-\mu=\sum_{i=1}^n r_i\alpha_i$\/ for some
non-negative integers $r_i$.

The elements $\lambda_i =\sum_{j=1}^i\epsilon_j$, $1\le i\le n$, are called
fundamental weights and the corresponding irreducible representations
$V(\lambda_i)$\/ the fundamental representations of \uqsl .

The representation  $V(\lambda_1)$\/ is called the natural representation of
\uqsl; it will be denoted by $V$ from now on. It has a basis $\{v_1, \ldots
,v_{n+1}\}$ on which the action is given by:
\begin{eqnarray*}
x_i^+.v_r&=& \delta_{r,i+1}v_{r-1},\\
x_i^-.v_r&=&\delta_{r,i}v_{r+1},\\
k_i.v_r&=&q^{\epsilon_r(i)} v_r\end{eqnarray*}
(we set $v_{-1}=v_{n+2}=0$).

Let $x_{\theta}^{{}\pm{}}$\/ be the operators on $V$\/ defined by
$$x_{\theta}^+.v_r =\delta_{r,n+1} v_1,\;\; x_{\theta}^-.v_r
=\delta_{r,1}v_{n+1},$$
and let $k_\theta=k_1k_2\ldots k_n$.
It is easy to see that $V$\/ can be made into a \uqsla--module $V(a)$, for all
$a\in\cit$, by letting $k_0$\/ act as $k_{\theta}^{-1}$ and
$x_0^{{}\pm{}}$\/ as $a^{{}\pm 1}x_{\theta}^{{}\mp{}}$.

\subsection{} \bed If $\ell\le n$, a finite--dimensional \uqsl--module $W$\/ is
 said to be of level $\ell$\/ if every irreducible component of $W$\/ is
isomorphic to an irreducible component of $V^{\otimes\ell}$.
\eed

Note that every level $\ell$\/ representation of \uqsl\/ is of type 1.

The next result follows immediately from Proposition 2.3 and the corresponding
classical result (which is well--known and easy to prove).

\bep Assume that $\ell\le n$. Then, the finite--dimensional \uqsl--module
$V(\lambda)$\/ is of level $\ell\le n$\/ iff $\sum_{i=1}^ni\lambda(i) =\ell$.
$\Box$
\eep
\brk This proposition shows that the concept of level is well--defined. The
assumption that $\ell\le n$\/ is necessary, for if $\ell_1$\/ or $\ell_2$\/ is
greater than $n$, it is possible for $V^{\ot \ell_1}$\/ and $V^{\ot \ell_2}$\/
to have an irreducible component in common even if $\ell_1\ne \ell_2$.
\erk
\subsection{}  It is easy to check that $c=k_0k_1\ldots k_n$\/ is central in
\uqsla .
\bep The central element $c$\/ of \uqsla\/  acts as $1$\/ on every
finite--dimensional \uqsla--module $W$\/ of type 1.
\eep
Proof. This was proved in [2] when $n=1$\/ and $W$\/ is irreducible.
Essentially the same proof works for all $n$\/ and the extension to arbitrary
finite--dimensional $W$\/ follows by an easy argument using Jordan--H\"older
series. $\Box$

\section{Hecke algebras and affine Hecke algebras}
In this section, we collect some well--known definitions and results concerning
(affine) Hecke algebras (cf. [9], [12]). We continue to assume that
$q\in\cit$\/ is not a root of unity.
\subsection{} \bed Fix $\ell\ge 1$. The affine Hecke algebra \hma\/ is the
unital associative algebra over $\ci$\/ with generators $\sigma_i^{{}\pm 1}$,
$i\in\{1,\ldots, \ell-1\}$, $y_j^{{}\pm 1}$, $j\in\{1,\ldots ,\ell\}$, and the
following defining relations:
\begin{eqnarray*}\sigma_i\sigma_i^{-1}&=&\sigma_i^{-1}\sigma_i =1,\\
\sigma_i\sigma_{i+1}\sigma_i &=&\sigma_{i+1}\sigma_i\sigma_{i+1},\\
\sigma_i\sigma_j&=&\sigma_j\sigma_i\;\;\;\; {\rm if}\;  |i-j| >1,\\
(\sigma_i +1)(\sigma_i -q^2) &=& 0,\\
y_jy_j^{-1}&=&y_j^{-1}y_j =1,\\
y_jy_k&=&y_ky_j,\\
y_j\sigma_i&=&\sigma_iy_j\;\;\;\; {\rm if}\; j\ne i\;{\rm or}\; i+1,\\
\sigma_iy_i\sigma_i &=& q^2y_{i+1} .\end{eqnarray*}

The unital associative algebra with generators $\sigma_i^{{}\pm 1}$,
$i\in\{1,\ldots ,\ell-1\}$, defined by the first four sets of relations above
is called the Hecke algebra \hm .
\eed

There is an obvious homomorphism of  \hm\/ onto the subalgebra of \hma\/
generated by the $\sigma_i$.
\bel  The multiplication map
 $C[y_1^{{}\pm 1},\ldots ,y_{\ell}^{{}\pm 1}]\ot \hm\to \hma$\/ is an
isomorphism of vector spaces.
$\Box$ \eel
\subsection{}  The following well-known result provides an analogue for affine
Hecke and Hecke algebras of the canonical homomorphism $S_{\ell_1}\times
S_{\ell_2}\to S_{\ell_1+\ell_2}$.
\bep There exists a unique homomorphism of algebras
$${\hat \iota}_{\ell_1,\ell_2}: {\hat H}_{\ell_1}(q^2)\otimes{\hat
H}_{\ell_2}(q^2)\to {\hat H}_{\ell_1+\ell_2}(q^2)$$
such that
\begin{eqnarray*}{\hat \iota}_{\ell_1,\ell_2}(\sigma_i\ot 1)= \sigma_i,
&\;\;&{\hat\iota}_{\ell_1,\ell_2}(y_j\ot 1)= y_j,
 \;\;i=1,\ldots,\ell_1-1,\;j=1,\ldots,\ell_1,\\
{\hat\iota}_{\ell_1,\ell_2}(1\ot \sigma_i)=\sigma_{i+\ell_1},&\;\;& {\hat
\iota}_{\ell_1,\ell_2}(1\ot y_j)= y_{j+\ell_1},
\;\;i=1,\ldots ,\ell_2-1,\;j=1,\ldots,\ell_2.\;\; \Box\end{eqnarray*}
\eep

Clearly the restriction of ${\hat\iota}_{\ell_1,\ell_2}$\/ to
$H_{\ell_1}(q^2)\otimes H_{\ell_2}(q^2)$\/ induces a homomorphism
$\iota_{\ell_1,\ell_2}:H_{\ell_1}(q^2)\otimes H_{\ell_2}(q^2)\to
H_{\ell_1+\ell_2}(q^2)$.

Let $M_i$\/ be a right $H_{\ell_i}(q^2)$--module for $i=1,2$, and let
$M_1\ot M_2$\/ be their outer tensor product (an $H_{\ell_1}(q^2)\otimes
H_{\ell_2}(q^2)$--module). Then, the $H_{\ell_1+\ell_2}(q^2)$--module $M_1\odot
M_2$, sometimes called the Zelevinsky tensor product of $M_1$\/ and $M_2$, is
defined by
$$M_1\odot M_2 = {\rm ind}_{H_{\ell_1}(q^2)\otimes
H_{\ell_2}(q^2)}^{H_{\ell_1+\ell_2}(q^2)}(M_1\ot M_2) = (M_1\ot
M_2)\bigotimes_{H_{\ell_1}(q^2)\otimes H_{\ell_2}(q^2)}
H_{\ell_1+\ell_2}(q^2).$$
The Zelevinsky tensor product $\hat{\odot}$\/ for affine Hecke algebra modules
is defined similarly. Standard properties of induced modules show that the
Zelevinsky tensor products are  associative up to isomorphism.

\subsection{} \bep Let $M_i$\/ be a finite--dimensional  $\hmai$--module,
$i=1,2$. Then, there is a canonical isomorphism of
$H_{\ell_1+\ell_2}(q^2)$--modules
$$(M_1\hat{\odot}M_2)\vert_{H_{\ell_1+\ell_2}(q^2)}\cong
M_1\vert_{H_{\ell_1}(q^2)}\odot M_2\vert_{H_{\ell_2}(q^2)},$$
where $M_i\vert_{H_{\ell_i}(q^2)}$\/ means $M_i$\/ regarded as an \hmi--module
by restriction, etc.
\eep
Proof. It is easy to see that the canonical map
$$M_1\vert_{H_{\ell_1}(q^2)}\odot
M_2\vert_{H_{\ell_2}(q^2)}\to(M_1\hat{\odot}M_2)\vert_{H_{\ell_1+\ell_2}(q^2)}$$
given by
$$(m_1\ot m_2)\ot h\mapsto (m_1\ot m_2)\ot h\;\;\;(m_i\in M_i,\;h\in
H_{\ell_1+\ell_2}(q^2))$$
is a well--defined surjective homomorphism of
$H_{\ell_1+\ell_2}(q^2)$--modules. But, by Lemma 3.1, the rank of
$\hat{H}_{\ell_1+\ell_2}(q^2)$\/ as an
$\hat{H}_{\ell_1}(q^2)\ot\hat{H}_{\ell_2}(q^2)$--module is the same as that of
$H_{\ell_1+\ell_2}(q^2)$\/ as an $H_{\ell_1}(q^2)\ot H_{\ell_2}(q^2)$--module.
It follows that
$$\dim_{\ci}(M_1\hat\odot M_2)=\dim_{\ci}(M_1\odot M_2).\;\;\;\Box
$$
\subsection{} Affine Hecke algebras have a family of universal modules, defined
as follows. Let $\una =(a_1,a_2,\ldots ,a_{\ell})\in(\cit)^{\ell}$\/ and set
$$M_{{\bf a}} =\hma/H_{{\bf a}} ,$$
the quotient of \hma\/ by the right ideal $H_{{\bf a}}$\/ generated by $y_j
-a_j$, $j=1,\ldots ,\ell$.
\bep {\rm ([12])}

\noindent (a) Every finite--dimensional irreducible \hma--module is isomorphic
to a quotient of some $M_{{\bf a}}$.

\noindent (b) For all ${\bf a}\in(\cit)^{\ell}$, $M_{{\bf a}}$\/ is isomorphic
as an \hm--module to the right regular representation.

\noindent (c) $M_{{\bf a}}$\/ is reducible as an \hma--module iff $a_j
=q^2a_k$\/ for some $j$, $k$. $\Box$\eep

\section{ Duality between \uqsla\/ and \hma}
We begin by recalling the duality, established by Jimbo [7], between
representations of \uqsl\/ and \hm .

\subsection{} Let $V$\/ be the natural $(n+1)$--dimensional representation of
\uqsl\/ defined in 2.4, and let $\check{R}:V\ot V\to V\ot V$\/ be the linear
map given by
\begin{equation}\check{R}(v_r\ot v_s) =\left\{\begin{array}{ll}
q^2 v_r\ot v_s& {\rm if}\; r =s,\\
qv_s\ot v_r& {\rm if}\; s >r,\\
qv_s\ot v_r + (q^2-1)v_r\ot v_s & {\rm if}\;
r>s.\end{array}\right.\end{equation}
 Fix $\ell >1$\/ and let $\check{R}_i\in {\rm End}_{\ci}(V^{\ot \ell})$\/ be
the map which acts as $\check{R}$\/ on the $i^{th}$\/ and $(i+1)^{th}$\/
factors of the tensor product, and as the identity on the other factors.

\bep{{\rm ([7])}} Fix $\ell,n\ge 1$. There  is a unique left \hm--module
structure on \vm\/ such that $\sigma_i$\/ acts as $\check{R}_i$\/ for
$i=1,\ldots ,\ell-1$. Moreover, the action of \hm\/ commutes with the natural
action of \uqsl\/ on \vm .

If $M$\/ is a right \hm--module, define
$${\cal J}(M) = M\otimes_{H_\ell(q^2)}\vm,$$
equipped with the natural left \uqsl--module structure induced by that on \vm .
Then, if $\ell\le n$, the functor $M\to {\cal J}(M)$\/ is an equivalence from
the category of finite--dimensional \hm--modules to the category of
finite-dimensional \uqsl--modules of level $\ell$. $\Box$
\eep

\subsection{} We now state the main result of this section, which is an
analogue of Proposition 4.1 for quantum affine algebras. Recall the operators
$k_{\theta}$, $x_{\theta}^{{}\pm{}}\in{\rm End}_{\ci}(V)$\/ defined in Section
2.4.
\bet Fix $\ell,n\ge 1$. There is a functor ${\cal F}$\/ from the category of
finite--dimensional right \hma--modules to the category of finite--dimensional
left \uqsla--modules of type 1 which are
of level $\ell$\/ as \uqsl--modules, defined as follows. If $M$\/ is an
\hma--module, then
${\cal F}(M) ={\cal J}(M)$\/ as a \uqsl
--module and the action of the remaining generators of \uqsla\/ is given by
\begin{eqnarray} x_0^{{}\pm{}}.(m\ot {\bf v}) &=&\sum_{j=1}^\ell m.y_j^{{}\pm
1}\ot Y_j^{{}\pm{}}.{\bf v},\\
k_0.(m\ot{\bf v}) &=&m\ot (k_{\theta}^{-1})^{\ot \ell}.{\bf v},\end{eqnarray}
where $m\in M$, ${\bf v}\in\vm$\/ and the operators $Y_j^{{}\pm{}}\in {\rm
End}_{\ci}\;(\vm)$, $j=1,\ldots,\ell$, are defined by
\begin{eqnarray*} Y_j^+&=&1^{\ot j-1}\ot x_{\theta}^-\ot (k_{\theta}^{-1})^{\ot
\ell -j},\\
Y_j^-&=& k_{\theta}^{\ot j-1}\ot x_{\theta}^+\ot 1^{\ot \ell
-j}.\end{eqnarray*}

The functor ${\cal F}$\/ is an equivalence of categories if $\ell\le n$.

\eet
Proof. We first show that the formulas (2) and (3) are well--defined. We do
this for the action of $x_0^+$, leaving the verification for $x_0^-$ and
$k_0$\/ to the reader.
Thus, we must prove that
$$x_0^+.(m.\sigma_i\ot{\bf v}) = x_0^+.(m\ot \sigma_i.{\bf v})$$
for $i=1,\ldots ,\ell$, ${\bf v}\in\vm$. This is equivalent to proving that, as
operators on
${\cal J}(M) = M\otimes_{H_\ell(q^2)}\vm$,
\begin{equation} \sum_{j=1}^{\ell}\sigma_iy_j\ot Y_j^+ = \sum_{j=1}^\ell
y_j\otimes Y_j^+\sigma_i .\end{equation}
If $j\ne i,i+1$, the $j^{th}$\/ terms on the left and right-hand sides of (4)
are equal, since $\sigma_iy_j = y_j\sigma_i$\/ and $\sigma_iY_j^+ =
Y_j^+\sigma_i$. Hence we must show that
$$\sigma_iy_i\ot Y_i^+ +\sigma_iy_{i+1}\ot Y_{i+1}^+=y_i\ot
Y_i^+\sigma_i+y_{i+1}\ot Y_{i+1}^+\sigma_i.$$
Using the relation $\sigma_i-(q^2-1)=q^2\sigma_i^{-1}$, this reduces to
$$q^2y_{i+1}\ot(\sigma_i^{-1}Y_i^+-Y_{i+1}^+\sigma_i^{-1})+y_i\ot(\sigma_iY_{i+1}^+-Y_i^+\sigma_i)=0.$$
Thus, it suffices to prove that
$$\sigma_iY_{i+1}^+ = Y_i^+\sigma_i,$$
i.e. that
\begin{equation}\check{R}(1\ot x_{\theta}^-) = (x_{\theta}^-\ot
k_{\theta}^{-1})\check{R}\end{equation}
as operators on $V\ot V$.
But this is easily checked by using the formula for $\check{R}$\/ in (1) and
that for $x_{\theta}^-$ in 2.4.

In proving that the formulas (2) and (3) define a representation of \uqsla, we
shall assume that $n >1$. The proof for the $sl_2$\/ case is similar (the
difference arises because the Dynkin diagram of $\hat{sl}_2$\/ has a double
bond).

The only relations to be checked are those involving $x_0^+$, $x_0^-$ and
$k_0$. This is straightforward except for the quantized Serre relations:
\begin{equation}
{}[x_i^{{}\pm{}},[x_0^{{}\pm{}}, x_i^{{}\pm{}}]_{q^{1/2}}]_{q^{1/2}} =
0,\end{equation}
\begin{equation}{}[x_0^{{}\pm{}},[x_i^{{}\pm{}},
x_0^{{}\pm{}}]_{q^{1/2}}]_{q^{1/2}} = 0,\end{equation}
for $i=1,n$. We verify (7) for $x_1^+$, leaving the other cases to the reader.

Applying the left--hand side of (7) to ${\cal J}(M)$\/ and considering the
terms involving $y_jy_k$, one sees that it suffices to prove that
\begin{equation} [Y_j^+,[\Delta^{(\ell)}(x_1^+), Y_k^+]_{q^{1/2}}]_{q^{1/2}} +
(j\leftrightarrow k) = 0,\end{equation}
where $(j\leftrightarrow k)$\/ means the result of interchanging $j$\/ and
$k$\/ in the first term and $\Delta^{(\ell)}$\/ is the $\ell^{th}$ iterated
comultiplication (so that $\Delta^{(2)}=\Delta$).
Equation (8) will be proved by induction on $\ell$, and we accordingly denote
$Y_k^+$\/ by $Y_k^{+(\ell)}$. If $\ell=1$, then (8) becomes
$$[x_{\theta}^{{}-{}},[x_1^{{}-{}}, x_{\theta}^{{}-{}}]_{q^{1/2}}]_{q^{1/2}} =
0,$$
which holds by the remarks at the end of 2.4.

For the inductive step we distinguish three cases:

\noindent (i) $j,k < \ell$,

\noindent (ii) $j <\ell$, $k =\ell$\/ or $j=\ell$, $k<\ell$,

\noindent (iii) $j=k=\ell$.

For the first case, notice that the left--hand side of (8) is
\begin{eqnarray*}[Y_j^{+(\ell-1)}\ot k_{\theta}^{-1},
[\Delta^{(\ell-1)}(x_1^+)\ot k_1 +1\ot x_1^+, Y_k^{+(\ell-1)}\ot
k_{\theta}^{-1}]_{q^{1/2}}]_{q^{1/2}} +(j\leftrightarrow k)&&\\
 =[Y_j^{+(\ell-1)}, [\Delta^{(\ell-1)}(x_1^+) +1\ot x_1^+,
Y_k^{+(\ell-1)}]_{q^{1/2}}]_{q^{1/2}}\ot k_1k_{\theta}^{-2} +(j\leftrightarrow
k)&&\\
+[Y_j^{+(\ell-1)}\ot k_{\theta}^{-1},  Y_k^{+(\ell-1)}\ot
[x_1^+,k_{\theta}^{-1}]_{q^{1/2}}]_{q^{1/2}} +(j\leftrightarrow
k).&&\end{eqnarray*}
The sum of the first two terms on the right-hand side vanishes by the induction
hypothesis, and the sum of the last two terms is a multiple of
$$[Y_j^{+(\ell-1)}\ot k_{\theta}^{-1},  Y_k^{+(\ell-1)}\ot
k_{\theta}^{-1}x_1^+]_{q^{1/2}}
 +(j\leftrightarrow k)=
 q^{1/2}[Y_j^{+(\ell-1)},Y_k^{+(\ell-1)}]\ot k_{\theta}^{-2}x_1^+
+(j\leftrightarrow k).$$
But the expression on the right-hand side is zero since
$[Y_j^{+(l-1)},Y_k^{+(l-1)}] =0$, so the induction step is established in this
case. The other two cases are similar; we omit the details.

We have thus proved that formulas (2) and (3) define a representation of
\uqsla. If $f:M\to M'$\/ is a homomorphism of \hma--modules, we define
${\cal F}(f):{\cal F}(M)\to {\cal F}(M')$\/ by
$${\cal F}(f)(m\ot{\bf v}) = f(m)\ot {\bf v}.$$
The proof that ${\cal F}(f)$\/ is a well--defined homomorphism of
\uqsla--modules is completely straightforward. It is now obvious that ${\cal
F}$\/ is a functor between the appropriate categories of representations.

\subsection{} Assume for the remainder of the proof that $\ell\le n$. To prove
that ${\cal F}$\/ is an equivalence, we must prove that

\noindent (a)  every finite--dimensional \uqsla--module $W$\/ of type 1 which
is of level $\ell$\/ as a \uqsl--module is isomorphic to ${\cal F}(M)$\/ for
some \hma--module $M$;

\noindent  (b) ${\cal F}$\/ is bijective on sets of morphisms.

\noindent (See [11], p.91.)

To prove (a), note that by Proposition 4.1, we may assume that  $W ={\cal
J}(M)$\/ for some \hm--module $M$. We shall reconstruct the action of the
$y_j^{{}\pm 1}$\/ on $M$\/ from the known action of $x_0^{{}\pm{}}$\/ and
$k_0$\/ on $W$.

We need the following lemma.

\bel (a) Let $M$\/ be a finite--dimensional \hm--module, and let ${\bf
v}\in\vm$. The linear map $M\to{\cal J}(M)$\/ given by $m\to m\ot{\bf v}$\/ is
injective if ${\bf v}$\/ has non--zero component in each isotypical component
of ${\cal J}(M)$.

\noindent (b) If $\{v_1,\dots,v_{n+1}\}$\/ is the standard basis of $V$,
$i_1,\ldots,i_\ell\in\{1,\ldots,n+1\}$\/ are distinct, and ${\bf
v}=v_{i_1}\ot\cdots\ot v_{i_\ell}$, then $V^{\ot\ell}=\uqsl.{\bf v}$. In
particular, \unv\/ satisfies the condition in part (a).
\eel
Proof. Part (a) follows easily from Proposition 4.1, and part (b) is
elementary. $\Box$

\subsection{} For $1\le j\le n$, let
\begin{eqnarray*}{\bf v}^{(j)}& =&v_2\ot\cdots\ot v_j\ot v_{n+1}\ot
v_{j+1}\ot\cdots \ot v_{\ell},\\
{\bf w}^{(j)}& =&v_2\ot\cdots\ot v_j\ot v_1\ot v_{j+1}\ot\cdots \ot
v_{\ell}.\end{eqnarray*}
Let ${\bf w}_{\tau}^{(j)}$\/ be the result of permuting the factors of ${\bf
w}^{(j)}$ by $\tau\in\sg$. Since $\{{\bf w}_{\tau}^{(j)}\}_{\tau\in\sg}$\/
clearly spans the subspace of \vm\/ of weight $\lambda_{\ell}$, we get, for any
$m\in M$,
$$x_0^-.(m\ot{\bf v}^{(j)} )=\sum_{\tau\in\sg}m_{\tau}\ot{\bf w}_{\tau}^{(j)}$$
for some $m_{\tau}\in M$. By (1), ${\bf w}_{\tau}^{(j)}$\/ is a (non--zero)
scalar multiple of $\sigma.{\bf w}^{(j)}$\/ for some $\sigma\in\hm$\/
(depending on $\tau$). It follows that
$$x_0^-.(m\ot{\bf v}^{(j)})=m'\ot{\bf w}^{(j)}$$
for some $m'\in M$. By Lemma 4.3, there exists $\alpha_j^-\in{\rm
End}_{\ci}(M)$\/ such that $m'=\alpha_j^-(m)$\/ for all $m\in M$. By  a similar
argument, there exists $\alpha_j^+\in{\rm End}_{\ci}(M)$\/ such that
\begin{eqnarray*}
x_0^+.(m\ot v_{n-\ell+2}\ot &\cdots &\ot v_{n-\ell+j}\ot v_1\ot
v_{n-\ell+j+1}\ot\cdots\ot v_n)\\
&=&\alpha_j^+(m)\ot v_{n-\ell+2}\ot\cdots\ot v_{n-\ell+j}\ot v_{n+1}\ot
v_{n-\ell+j+1}\ot\cdots\ot v_n \end{eqnarray*}
for all $m\in M$.
\subsection{} We need to prove the following lemma. The proof of the theorem
itself continues in Section 4.6.
\bel For all $m\in M$, ${\bf v}\in\vm$, we have
$$x_0^{{}\pm{}}.(m\ot{\bf v}) =\sum_{j=1}^\ell\alpha_j^{{}\pm{}}(m)\ot
Y_j^{{}\pm{}}.{\bf v}.$$
\eel
Proof. Let $\unv =v_{i_1}\ot\cdots\ot v_{i_{\ell}}$. If $\{i_1,\ldots
,i_{\ell}\}\subset \{1,\ldots ,n\}$, it is clear that $x_0^-.(m\ot\unv) = 0$,
since $\epsilon_{i_1}+\ldots +\epsilon_{i_{\ell}}+\epsilon_1+\ldots
+\epsilon_n$\/ cannot be a weight of \vm .

Let $r\ge 0$, $s\ge 1$, $1\le j_1<j_2<\ldots <j_r\le\ell$, $1\le j'{}_1<
j'{}_2<\ldots <\
j'{}_s\le\ell$, and assume that
$\{j_1,\ldots ,j_{\ell}\}\cap\{j'{}_1,\dots ,j'{}_s\} =\emptyset$. Write ${\bf j}
=(j_1,\dots , j_r)$, ${\bf j}' =(j'{}_1,\dots ,j'{}_s)$, and let $V^{({\bf j},{\bf
j}')}$\/ be the subspace of \vm\/ spanned by vectors which have $v_1$\/ in
positions $j_1,\ldots, j_r$, $v_{n+1}$\/ in positions $j'{}_1,\ldots, j'{}_s$, and
vectors from $\{v_2,\ldots, v_n\}$\/ in the remaining positions. We shall prove
the lemma when ${\bf v}\in V^{({\bf j},{\bf j}')}$\/ for all such ${\bf j},{\bf
j}'$\/ in two steps:

\noindent (i) for $s=1$, by induction on $r$;

\noindent (ii) for all $r$, by induction on $s$.

Observe that, by Lemma 4.3 (b) applied to the subalgebra of \uqsl\/ generated
by the $x_i^{{}\pm{}}$, $k_i^{{}\pm 1}$\/ for $i\in\{2,\ldots ,n\}$, to prove
Lemma 4.5 for all ${\bf v}\in V^{({\bf j},{\bf j}')}$, it suffices to prove it
for one ${\bf v}\in V^{({\bf j},{\bf j}')}$\/ with the property that no vector
from the set $\{v_2,\ldots ,v_n\}$\/ is repeated. (Note that such vectors ${\bf
v}$\/ exist since $\ell+1-r-s\le\ell\le n$.)
\vskip 6pt
Proof of Step (i). If $r =0$ (and $s=1$), there is nothing to prove, for we can
take $\unv =v_2\ot\cdots\ot v_{j'{}_1}\ot v_{n+1}\ot v_{j'{}_1+1}\ot\cdots\ot
v_{\ell}$\/ and use the definition of $\alpha_{j'{}_1}^-$. Assume that the result
holds for $r-1$, and let ${\tilde\unj}  =(j_1,\ldots ,j_{r-1})$. Let ${\bf
v}'\in V^{(\tilde{\bf j}, j')}$\/ have $v_2$\/ in the $j_r^{th}$\/ position,
and distinct vectors from $\{v_3,\ldots ,v_n\}$\/ in the remaining positions.
Then,
$$\unv = x_1^+.\unv' .$$ Let $\unv''$\/ (resp. $\unv'''$) be the element
obtained from $\unv'$\/ by replacing $v_{n+1}$\/ by $v_1$\/ (resp. $v_2$\/ by
$v_1$).
We then get, for all $m\in M$,
\begin{eqnarray*}x_0^-.(m\ot\unv) &= &x_1^+x_0^-.(m\ot {\bf v}')\\
&=& q^{|\{t<r|j_t<j'{}_1\}|}\alpha^-_{j'{}_1}(m)\ot (1^{\ot j_{r}-1}\ot x_1^+\ot
k_1^{\ell -j_r }).\unv''\\
&=&q^{|\{t<r|j_t<j'{}_1\}|}q^{\delta_{j_r<j_1'}}\alpha^-_{j'{}_1}(m)\ot \unv'''\\
&=&q^{|\{t\le r|j_t<j'{}_1\}|}\alpha_{j'{}_1}^-(m)\ot \unv'''\\
&=&\alpha_{j'{}_1}^-(m)\ot Y_{j'{}_1}^-.\unv .\end{eqnarray*}

Proof of Step (ii). Assume that the result holds for all ${\bf v}\in V^{({\bf
j},{\bf j}')}$\/ with fewer than $s$\/ $v_{n+1}$s. It suffices, as in step 1,
to prove the result for one element $\unv\in V^{({\bf j} ,{\bf j}')}$\/ which
has distinct entries from $\{v_3,\ldots ,v_n\}$
in the remaining positions. Fix such a $\unv$\/ and let $\unv'$\/ be the
element obtained from \unv\/ by replacing $v_{n+1}$\/ in positions $j_1'$\/ and
$j_2'$\/ by $v_n$.
Then,
$$\unv =\frac{(x_n^-)^2}{q+q^{-1}}.\unv' .$$
Using a quantized Serre relation we get
$$ x_0^-.(m\ot \unv) = x_n^-x_0^-x_n^-.(m\ot\unv')
-\frac{(x_n^-)^2x_0^-}{q+q^{-1}}.(m\ot\unv').$$
Since $x_n^-$\/ operates in the $j'{}_{1}^{th}$\/ and $j'{}_2^{th}$\/ positions in
$\unv'$, we obtain, using the induction hypothesis,
$$\frac{(x_n^-)^2x_0^-}{q+q^{-1}}.(m\ot\unv ') =q^2\sum_{k=3}^s
\alpha_{j_k'}^-(m)\ot Y_{j'{}_k}^-.\unv .$$
On the other hand,
$$x_n^-.(m\ot \unv') =m\ot\unv'' +q^{-1}m\ot \unv''',$$
where $\unv''$\/ (resp. $\unv'''$) is obtained from $\unv'$\/ by replacing
the $v_n$\/ in its $j'{}_1^{th}$\/ position (resp. $j'{}_2^{th}$\/ position) by
$v_{n+1}$. Using the induction hypothesis, we get
$$x_0^-x_n^-.(m\ot\unv') =\sum_{k\ne 2}\alpha_{j'{}_k}^-(m)\ot Y_{j'{}_k}^-.\unv''
+q^{-1}\sum_{k\ne 1}\alpha_{j'{}_k}^-(m)\ot Y_{j'{}_k}^-.\unv'''.$$
Noting that $\unv'''$\/ has $v_n$\/ only in the $j'{}_2^{th}$\/ position, we find
that
$$x_n^-.\sum_{k\ne 2}\alpha_{j'{}_k}^-(m)\ot Y_{j'{}_k}^-.\unv'' =
\alpha_{j'{}_1}^-(m)\ot Y_{j'{}_1}^-.\unv' + q^2\sum_{k >2} \alpha_{j'{}_k}^-(m)\ot
Y_{j'{}_k}^-.\unv'.$$
Similarly,
$$x_n^-.\sum_{k\ne 1}\alpha_{j'{}_k}^-(m)\ot Y_{j'{}_k}^-.\unv''' =   q\sum_{k \ne
1} \alpha_{j'{}_k}^-(m)\ot Y_{j'{}_k}^-.\unv'.$$
Combining these computations we obtain finally,
\begin{eqnarray*} x_0^-.(m\ot\unv)&=&-q^2\sum_{k >2} \alpha_{j'{}_k}^-(m)\ot
Y_{j'{}_k}^-.\unv\\
&+&q^2\sum_{k >2} \alpha_{j'{}_k}^-(m)\ot Y_{j'{}_k}^-.\unv + \alpha_{j'{}_1}^-(m)\ot
Y_{j'{}_1}^-.\unv\\
&+&\sum_{k\ne 1} \alpha_{j'{}_k}^-(m)\ot Y_{j'{}_k}^-.\unv\\
&=&\sum_{k =1}^s \alpha_{j'{}_k}^-(m)\ot Y_{j'{}_k}^-.\unv ,\end{eqnarray*}
as required.

This proves Lemma 4.5 for $x_0^-$. The proof for $x_0^+$\/ is similar. $\Box$

\subsection{} We can now complete the proof of the theorem. We show that
setting $$m.y_j^{{}\pm 1} =\alpha_j^{{}\pm{}}(m)$$ defines
a right \hma--module structure on $M$, extending its \hm--module structure. We
have to check the following relations:

\noindent (i) $y_jy_j^{-1} =y_j^{-1}y_j =1$,

\noindent (ii) $y_jy_k =y_ky_j$,

\noindent (iii) $q^2y_{j+1} =\sigma_jy_j\sigma_j$.

Relations (i) and (ii) are proved by computing both sides of  the equation
$$[x_0^+,x_0^-].(m\ot{\bf v}) =\left(\frac{k_0-k_0^{-1}}{q-q^{-1}}\right).(m\ot
{\bf v}),$$
where in the first case we take ${\bf v}$\/ to be a vector with $v_{n+1}$\/ in
the $j^{th}$\/ place and $v_{n-\ell+2},\ldots ,v_{n}$\/  in the remaining
places (in any order), and in the second case we take ${\bf v}$\/ to be a
vector with $v_1$\/ in the $j^{th}$\/ place, $v_{n+1}$\/ in the $k^{th}$\/
place and distinct vectors from $\{v_2,\ldots ,v_n\}$\/ in the other places.
Notice that since the central element $c\in\uqsla$\/ acts as 1 on $W$ we have
$k_0.(m\ot{\bf v}) =m\ot(k_{\theta}^{-1})^{\ot\ell}.{\bf v}$.

To prove (iii), let $\unv = v_{i_1}\ot\cdots\ot v_{i_{\ell}}\in\vm$, where $i_j
=2$, $i_{j+1} =1$, and the remaining $i_k$\/ are distinct elements from
$\{3,\ldots ,n\}$\/ (this is possible since $\ell\le n$). Let $\unv '$\/ be the
result of replacing $v_1$\/ in the $i_{j+1}^{th}$\/ position in \unv\/ by
$v_{n+1}$. Since
$$\check{R}(v_2\ot v_{n+1}) = qv_{n+1}\ot v_2,\;\; \check R(v_1\ot v_2) =
qv_2\ot v_1,$$
we have, for all $m\in M$,
$$m.\sigma_jy_j\sigma_j\ot\unv ' = qm.\sigma_jy_j\ot\unv '',$$
where $\unv ''$\/ is obtained from $\unv'$\/ by interchanging its $j^{th}$\/
and $(j+1)^{th}$\/ factors, which
$$ = qx_0^+.(m.\sigma_j\ot\unv '''),$$
where $\unv '''$\/ is obtained from \unv\/ by interchanging its $j^{th}$\/ and
$(j+1)^{th}$\/ factors, which
$$ =q^2x_0^+.(m\ot \unv) = q^2m.y_{j+1}\ot\unv '.$$
Since $\unv '$\/ has distinct components, Lemma 4.3 implies that
$$q^2m.y_{j+1} = m.\sigma_jy_j\sigma_j,$$
for all $m\in M$.

The proof that $W\cong {\cal F}(M)$\/ as \uqsla--modules is now complete. To
show that ${\cal F}$\/ is an equivalence, we must prove that it is bijective on
sets of morphisms. Injectivity of ${\cal F}$\/ follows from that of ${\cal J}$.
For surjectivity, let $F:{\cal F}(M)\to{\cal F}(M')$\/ be a homomorphism of
\uqsla--modules. By Proposition 4.1 again, $F ={\cal J}(f)$\/ for some
homomorphism $f:M\to M'$\/ of \hm--modules. The fact that $F$\/ commutes with
the action of $x_0^+$\/ gives
$$\sum_{j=1}^{\ell}f(m.y_j)\ot Y_j^+.\unv =\sum_{j=1}^{\ell}f(m).y_j\ot
Y_j^+.\unv$$
for all $m\in M$, $\unv\in\vm$. By choosing \unv\/ suitably, as in the
preceding part of the proof, it is easy to see that this implies
$$f(m.y_j) = f(m).y_j$$
for all $j=1,\ldots ,\ell$. $\Box$
\subsection{} The functor ${\cal F}$\/ is clearly one of \ci--linear
categories. The following result shows that it also captures part of the tensor
structure of the category of \uqsla--modules.
\bep Let $M_i$\/ be a finite--dimensional \hmai--module, $i=1,2$. Then, there
is a canonical isomorphsm of \uqsla--modules
$${\cal F}(M_1\hat{\odot} M_2)\cong {\cal F}(M_1)\ot {\cal F}(M_2).$$
\eep
Proof. We recall the following elementary fact: if $\iota: B\to A$\/ is a
homomorphism of unital associative algebras over a field, $M$\/ is a right
$B$--module, $W$\/ a left $A$--module, and $W\vert_B$\/ is $W$\/ regarded as a
left $B$--module via $\iota$, there is a canonical isomorphism of vector spaces
$${\rm ind}_B^A(M)\ot W\cong M\bigotimes_B W\vert_B.$$
In fact, the isomorphism is given by
$$(m\ot a)\ot w\to m\ot aw \;\;\;\;\;(m\in M,a\in A, w\in W).$$

Taking $A= H_{\ell_1+\ell_2}(q^2)$, $B =\hmo\ot\hmt$,
$\iota=\iota_{\ell_1,\ell_2}$, $M =M_1\ot M_2$\/ and $W =
V^{\ot\ell_1+\ell_2}$, and noting that
$W\cong (V^{\ot\ell_1})\ot (V^{\ot\ell_2})$\/ as an $\hmo\ot\hmt$--module, we
get a canonical isomorphism of vector spaces
$${\cal F}(M_1\hat\odot M_2)\to (M_1\ot
M_2)\bigotimes_{\hmo\ot\hmt}(V^{\ot\ell_1}\ot V^{\ot\ell_2}).$$
The right--hand side is obviously isomorphic to ${\cal F}(M_1)\ot {\cal
F}(M_2)$\/ as a vector space. To complete the proof, one must check that the
resulting isomorphism of vector spaces
$${\cal F}(M_1\hat\odot M_2)\to {\cal F}(M_1)\ot {\cal F}(M_2)$$
commutes with the action of \uqsla . This is completely straightforward. $\Box$
\subsection{} We analyze the functor ${\cal F}$\/ in more detail in Section 7,
when the parametrizations of the finite--dimensional irreducible
representations of \hma\/ and \uqsla\/ have been described. The following
result is, however, easy to prove now.. Recall the universal \hma--modules
$M_{{\bf a}}$\/ and the \uqsla--modules $V(a)$\/ defined in Sections 2.4 and
3.4 respectively.
\bep Let $\una =(a_1,\ldots ,a_\ell)\in(\cit)^\ell$, $\ell,n\ge 1$. There is a
canonical isomorphism of \uqsla--modules
$${\cal F}(M_{{\bf a}})\cong V(a_1)\ot\cdots\ot V(a_{\ell}).$$
\eep
Proof. As an \hm--module, $M_{{\bf a}}$\/ is the right regular representation.
It follows that the map
\begin{equation} \vm\to {\cal J}(M_{{\bf a}})\end{equation}
given by $\unv\to 1\ot\unv$\/ is an isomorphism of \uqsl--modules. Now,
$$x_0^+.(1\ot\unv) =\sum_{j=1}^{\ell} 1.y_j\ot Y_j^+.\unv =
(\sum_{j=1}^{\ell}a_jY_j^+).\unv .$$
On the other hand,
$$\Delta^{(\ell)}(x_0^+) =\sum_{j=1}^{\ell}1^{\ot j-1}\ot x_0^+\ot k_0^{\ot
\ell-j}$$
acts on $V(a_1)\ot\cdots\ot V(a_{\ell})$\/ as
$$\sum_{j=1}^{\ell}1^{\ot j-1}\ot a_jx_{\theta}^-\ot (k_{\theta}^{-1})^{\ot
\ell-j} = \sum_{j=1}^{\ell}a_jY_j^+.$$
One checks in the same way that the map in (9) commutes with the action of
$x_0^-$\/ and $k_0$. $\Box$
\becor Let $1\le \ell\le n$.

\noindent (a) Every finite--dimensional \uqsla--module of type 1 and level
$\ell$\/ as a \uqsl--module is isomorphic to a quotient of $V(a_1)\ot\cdots\ot
V(a_{\ell})$, for some $a_1,\ldots ,a_{\ell}\in\cit$.

\noindent (b) If $a_1,\ldots, a_{\ell}\in\cit$, then $V(a_1)\ot\cdots\ot
V(a_{\ell})$\/ is reducible as a \uqsla--module iff $a_j =q^2a_k$\/ for some
$j,k$.
\eecor
Proof. This  follows immediately from Proposition 3.4 and  the fact that ${\cal
F}$\/ is an equivalence of categories. $\Box$

\subsection{} Theorem 4.2 has a classical analogue, in which \uqsla\/ is
replaced by (the universal enveloping algebra of) the affine Lie algebra
${\hat{sl}_{n+1}}$, and \hma\/ by (the group algebra of) the affine Weyl group
of $GL_{\ell}(\ci)$, i.e. the semi--direct product $S_{\ell}{\tilde{\times}}
{\bf Z^{\ell}}$, where \sg\/ acts on the additive group ${\bf Z^{\ell}}$\/ by
permuting the coordinates. We recall that ${\hat {sl}}_{n+1}$\/ is the
universal central extension (with one--dimensional centre) of the Lie algebra
$L(sl_{n+1})$\/ of Laurent polynomial maps $\cit\to sl_{n+1}$. We identify
$sl_{n+1}$\/ with the subalgebra of $L(sl_{n+1})$\/ consisting of the constant
maps.
\bet There is a functor ${\cal F}_0$\/ from the category of finite--dimensional
 $\sg{\tilde{\times}} {\bf Z}^\ell$--modules to the category of
finite--dimensional $L(sl_{n+1})$--modules which are of level $\ell$\/ as
$sl_{n+1}$--modules, defined as follows. One takes
$${\cal F}_0(M) = M\bigotimes_{\sg}\vom$$
with the action of $f\in L(sl_{n+1})$\/ given by
$$f.(m\ot\unv) =\sum_{j=1}^{\ell} m.z_j\ot (1^{\ot{j-1 }}\ot f(1)\ot 1^{\ot
\ell-j}).\unv ,$$
where $z_j =(0,\ldots ,0,1,0,\ldots ,0)\in {\bf
Z^{\ell}}\subset\sg{\tilde{\times}} {\bf Z^{\ell}}$\/ (with $1$\/ in the
$j^{th}$\/ position).
If $\ell\le n$, ${\cal F}_0$\/ is an equivalence. $\Box$
\eet

The proof of this theorem is analogous to (but simpler than) that of Theorem
4.2.
\brk The finite--dimensional irreducible representations of $L(sl_{n+1})$\/
were classified in [1]. For any $a\in\cit$, there is a homomorphism of Lie
algebras
$$ev_a^0: L(sl_{n+1})\to sl_{n+1}$$
given by $ev_a^0(f) =f(a)$. If $W$\/ is an irreducible $sl_{n+1}$--module,
pulling back by $ev_a^0$\/ gives an irreducible $L(sl_{n+1})$--module $W(a)$.
It is not difficult to prove that every finite--dimensional irreducible
representation  of $L(sl_{n+1})$\/ is isomorphic to a tensor product of
$W(a)$s.

It is easy to identify the corresponding representations of
$\sg\tilde\times{{\bf Z}^{\ell}}$. There is a homomorphism
$$\tilde{ev}_a^0:\sg\tilde\times{{\bf Z^{\ell}}}\to\sg $$
which is the identity on \sg\/ and for which $\tilde{ev}_a^0(z_j) = a$\/ for
all $j$. If $M$\/ is an irreducible \sg--module, pulling $M$\/ back by
$\tilde{ev}_a^0$\/ gives an irreducible $\sg\tilde\times{{\bf
Z^{\ell}}}$--module $M(a)$. It is clear that
$${\cal F}_0(M(a))\cong {\cal FS}(M)(a).$$
By Theorem 4.9, every finite--dimensional irreducible $\sg\tilde\times{{\bf
Z^{\ell}}}$--module is isomorphic to a Zelevinsky tensor product of $M(a)$s.
\erk

\section{Evaluation Representations} In this section, we construct analogues
for \uqsla\/ and \hma\/ of the representations of $sl_{n+1}$\/ and
$S_{\ell}\tilde{\times}{\bf Z}^{\ell}$\/ described in Remark 4.9, and show how
these representations are related by the functor ${\cal F}$.

\subsection{} The following result was observed by Cherednik [4]. The proof is
straightforward.
\bep For every $a\in\cit$, there exists a  homomorphism
${\tilde{ev}}_a:\hma\to\hm$\/ such that
\begin{eqnarray*} \tilde{ev}_a(\sigma_i) &=&\sigma_i,\\
\tilde{ev}_a(y_j)&=& aq^{-2(j-1)}\sigma_{j-1}\sigma_{j-2}\ldots
\sigma_2\sigma_1^2\sigma_2\ldots \sigma_{j-1},\end{eqnarray*}
for $i=1,\ldots \ell-1$, $j=1,\ldots ,\ell$. $\Box$
\eep

Note that ${\tilde{ev}_a}$\/ can be characterized as the unique homomorphism
$\hma\to\hm$\/ which is the identity on $\hm\subset\hma$\/ and which maps
$y_1$\/ to $a$.

If $M$\/ is any \hm--module, pulling back $M$\/ by ${\tilde{ev}}_a$\/ gives an
\hma--module $M(a)$\/ which is isomorphic to $M$\/ as an \hm--module.
\subsection{} In [7], Jimbo defined a quantum analogue of the homomorphism
$ev_a^0:\hat{sl}_{n+1}\to sl_{n+1}$. To describe it, we need the following
\bed \uqgl\/ is the associative algebra over \ci\/ with generators
$x_i^{{}\pm{}}$, $i=1,\ldots ,n$, $t_r^{{}\pm 1}$, $r=1,\ldots ,n+1$, and the
following defining relations:
\begin{eqnarray*}
t_rt_r^{-1}=&1&=t_r^{-1}t_r,\\
t_rt_s&=&t_st_r\;,\\
t_rx_i^{{}\pm{}}t_r^{-1}&=& q^{{}\pm
(\delta_{r,i}-\delta_{r,i+1})}x_i^{{}\pm{}}\; ,\\
{}[x_i^{{}\pm{}},[x_j^{{}\pm{}},x_i^{{}\pm{}}]_{q^{1/2}}]_{q^{1/2}} &=&
0\;\;{\rm if}\;|i-j| =1,\\
{}[x_i^{{}\pm{}}, x_j^{{}\pm{}}]&=& 0\;\;{\rm if}\; |i-j|>1,\\
{}[x_i^+ , x_j^-]&=&\delta_{ij}\frac{k_i -k_i^{-1}}{q
-q^{-1}}\;,\end{eqnarray*}
where $k_i= t_it_{i+1}^{-1}$.
\eed

The algebra \uqgl\/ has a Hopf algebra structure, but we shall not make any use
of it.

Note that there is an obvious homomorphism $\uqsl\to\uqgl$.

\subsection{} Fix an $(n+1)^{th}$\/ root $q^{1/{(n+1)}}$\/ of $q$. We shall say
that  a finite--dimensional \uqgl--module $W$\/ is of type 1 if

\noindent (a) $W$\/ is of type 1 regarded as a \uqsl--module,

\noindent (b) the $t_r$\/ act semisimply on $W$\/ with eigenvalues which are
integer powers of $q^{1/{(n+1)}}$,

\noindent  (c) $t_1t_2\ldots t_{n+1}$\/ acts as 1 on $W$.

It is easy to see that restriction to \uqsl\/ is an equivalence from the
category of finite--dimensional \uqgl--modules of type 1 to the category of
finite--dimensional \uqsl--modules of type 1. In particular the functor ${\cal
J}$\/ of Proposition 4.1 may be viewed as taking values in the category of
finite--dimensional \uqgl--modules of type 1.
\subsection{} We can now state
\bep{\rm ([7])}
 For any $a\in\cit$, there exists a homomorphism $ev_a:\uqsla\to\uqgl$\/ such
that
\begin{eqnarray*}
ev_a(x_i^{{}\pm{}}) &=& x_i^{{}\pm{}},\;\; ev_a(k_i) = k_i,\;\; i=1,\ldots
,n,\\
ev_a(k_0)&=& (k_1k_2\ldots k_n)^{-1},\\
ev_a(x_0^{{}\pm{}})&=& ({{}\pm 1})^{(n-1)}q^{{}\mp {(n+1)/2}}a^{{}\pm
1}(t_1t_{n+1})^{{}\pm 1}[x_n^{{}\mp{}},[x_{n-1}^{{}\mp {}},\ldots
,[x_2^{{}\mp{}},x_1^{{}\mp{}}]_{q^{1/2}}\ldots ]_{q^{1/2}}]_{q^{1/2}}.\;\Box
\end{eqnarray*}
\eep

If $W$\/ is a \uqsl--module of type 1, we may regard $W$\/ as a \uqgl--module
by (5.3). The pull-back of $W$\/ by the homomorphism $ev_a$\/ is a
\uqsla--module which we denote by $W(a)$.

\subsection{} The main result of this section is
\bet Let $1\le \ell\le n$, and let $M$\/ be a finite--dimensional right
\hm--module. Then there is a  canonical isomorphism of \uqsla--modules,
$${\cal F}(M(q^{-{2\ell}/{(n+1)}}a))\cong {\cal J}(M)(a),$$
for all $a\in\cit$.
\eet
Proof. By Theorem 4.2 we know that ${\cal J}(M)(a)\cong {\cal F}(N)$, for some
\hma--module $N$\/ which is isomorphic to $M$\/ as an \hm--module. It suffices
to prove that $y_1$\/ acts as the scalar $a$\/ on $N$. To prove this, we
compute the action of $x_0^+$\/ on $m\ot v_1\ot v_{n-\ell+2}\ot
v_{n-\ell+3}\ot\cdots\ot v_n\in{\cal F}(N)$\/ in two different ways, for all
$m\in M$.

First, by the definition of ${\cal F}$, we have
\begin{equation}x_0^+.(m\ot v_1\ot v_{n-\ell+2}\ot v_{n-\ell+3}\ot\cdots\ot
v_n) =m.y_1\ot v_{n+1}\ot v_{n-\ell+2}\ot v_{n-\ell+3}\ot\cdots\ot v_n.
\end{equation}
On the other hand, let $f_n = [x_n^{{}-{}},[x_{n-1}^{{}- {}},\ldots
,[x_2^{{}-{}},x_1^{{}-{}}]_{q^{1/2}}\ldots ]_{q^{1/2}}]_{q^{1/2}}$.
Then,
\begin{eqnarray}
x_0^+.(m\ot v_1\ot v_{n-\ell+2}\ot\cdots \ot v_n)&=&
m\ot ev_a(x_0^+).(v_1\ot v_{n-\ell+2}\ot\cdots\ot v_n)\nonumber\\
&=& aq^{-(n-1)/2-2\ell/(n+1)} m\ot f_n.(v_1\ot v_{n-\ell+2}\ot\cdots\ot
v_n).\end{eqnarray}

We prove by induction on $n$\/ that
$$f_n.(v_1\ot v_{n-\ell+2}\ot\cdots\ot v_n) =q^{(n-1)/2} v_{n+1}\ot
v_{n-\ell+2}\ot\cdots\ot v_n.$$
The result is obvious if $n=1$. Assuming it for $n-1$, note that $f_n
=[x_n^-,f_{n-1}]_{q^{1/2}}$, so by the induction hypothesis,
\begin{eqnarray*} f_n.(v_1\ot v_{n-\ell+2}\ot\cdots\ot v_n)&=&
q^{(n-1)/2}x_n^-.(v_n\ot v_{n-\ell+2}\ot\cdots\ot v_n)\\
&&\quad -q^{-1/2}f_{n-1}.(v_1\ot v_{n-\ell+2}\ot\ldots\ot v_{n-1}\ot
v_{n+1}).\end{eqnarray*}
Since $x_i^-.v_{n+1} =0$\/ for $1\le i\le n-1$, we see that
\begin{eqnarray*}
f_{n-1}.((v_1\ot v_{n-\ell+2}\ot\cdots\ot v_{n-1})\ot v_{n+1})&=&
(f_{n-1}.(v_1\ot v_{n-\ell+2}\ot\cdots\ot v_{n-1}))\ot v_{n+1} \\
&=&q^{(n-2)/2}v_n\ot v_{n-\ell+2}\ot\cdots\ot v_{n-1}\ot
v_{n+1},\end{eqnarray*}
by the induction hypothesis again. Hence,
\begin{eqnarray*}
f_n.(v_1\ot v_{n-\ell+2}\ot\cdots\ot v_n)&=& q^{(n-1)/2}(v_{n+1}\ot
v_{n-\ell+2}\ot\cdots\ot v_n+q^{-1}v_n\ot v_{n-\ell+2}\ot\cdots\ot v_{n-1}\ot
v_{n+1})\\
&&-q^{(n-3)/2}v_n\ot v_{n-\ell+2}\ot\cdots\ot v_{n-1}\ot v_{n+1}\\
&=& q^{(n-1)/2}v_{n+1}\ot v_{n-\ell+2}\ot\cdots\ot v_n,
\end{eqnarray*}
as required.

Hence, from (11), we obtain
$$x_0^+.(m\ot v_1\ot v_{n-\ell+2}\ot\cdots\ot v_n)=aq^{-2\ell/(n+1)}m\ot
v_{n+1}\ot v_{n-\ell+2}\ot\cdots\ot v_n.$$
Comparing with (10), and using Lemma 4.3, we obtain
$$m.y_1=aq^{-2\ell/(n+1)}m$$
for all $m\in M$. $\Box$

\section{Classification of finite--dimensional \uqsla--modules}
\subsection{} The finite--dimensional irreducible \uqsla--modules of type 1
were classified in [2], [3]. To describe this result, we need an alternative
presentation of \uqsla\/ given in [6]. By Proposition 2.6, we need only
consider the quantum loop algebra \uqslp , the quotient of \uqsla\/ by the two
sided ideal generated by $c-1$.

\bep \uqslp\/ is isomorphic as an algebra to the algebra $\cal{A}$\/ with
generators $X_{i,r}^{{}\pm{}}$ ($ i\in\{1,\ldots, n\}, r\in {\bf Z}$),
$H_{i,r}$ ($i\in \{1,\ldots ,n\}$, $r\in {\bf Z}\backslash \{0\}$), and
$K_i^{{}\pm 1}$, ($i\in \{1,\ldots ,n\}$), and the following defining
relations:
\begin{eqnarray*}
K_iK_i^{- 1} =&1& =  K_i^{ -1}K_i,\\
K_iH_{j,r}&=&H_{j,r}K_i\;,\\
{}[H_{i,r}, H_{j,s}]&=& 0\;,\\
K_iX_j^{{}\pm{}}K_i^{-1}&=& q^{{}\pm a_{ij}}X_j^{{}\pm{}}\; ,\\
{}[H_{i,r}, X_{j,s}^{{}\pm{}}]&=&{{}\pm \frac1r}[ra_{ij}]_qX_{j,r+s}^{{}\pm
{}}\;,\\
X_{i,r+1}^{{}\pm{}}X_{j,s}^{{}\pm{}} - q^{{}\pm
a_{ij}}X_{j,s}^{{}\pm{}}X_{i,r+1}^{{}\pm{}}&=&q^{{}\pm a_{ij}}
X_{i,r}^{{}\pm{}}X_{j,s+1}^{{}\pm{}} - X_{j,s+1}^{{}\pm{}}X_{i,r}^{{}\pm{}},\\
{}[X_{i,r}^+ , X_{j,s}^-]&=&\delta_{ij}\frac{\Phi_{i,r+s}^+ -\Phi_{i,r+s}^-
}{q -q^{-1}}\;,\end{eqnarray*}
$$\sum_{\pi\in S_p}\sum_{k=0}^{p}
 (-1)^k\left[{p\atop k}\right]_{q}\; X_{i,r_{\pi(1)}}^{{}\pm{}}\ldots
X_{i,r_{\pi(k)}}^{{}\pm{}}X_{j,s}^{{}\pm{}}X_{i,r_{\pi(k+1)}}^{{}\pm{}}\ldots
X_{i,r_{\pi(p)}}^{{}\pm{}}= 0\;,\;i\ne j ,$$
for all sequences $(r_1,\ldots,r_p)\in{\bf Z}^p$,
where $p = 1-a_{ij}$\/ and the elements $\Phi_{i,r}^{{}\pm{}}$\/ are determined
by equating coefficients of powers of $u$\/ in the formal power series
$$\sum_{r=0}^{\infty}\Phi_{i, {{}\pm r}}^{{}\pm{}}u^{{}\pm r} = K_i^{{}\pm
1}exp({}\pm (q-q^{-1})\sum_{s=1}^{\infty} H_{i, {}\pm s} u^{{}\pm s}) .$$

The isomorphism $f:\uqslp\to{\cal A}$\/ is given by
$$f(x_i^{{}\pm{}}) = X_{i,0}^{{}\pm{}},\;\; f(k_i^{{}\pm 1}) = K_i^{{}\pm 1},$$
for $i\in\{1,\ldots, n\}$,
and
\begin{eqnarray*} f(k_0^{{}\pm 1}) &=& (K_1K_2\ldots K_n)^{{}\mp 1}, \\
f(x_0^+)&=& (-1)^{m-1}q^{-(n-3)/2}[X_{n,0}^-,[X_{n-1, 0}^-,\ldots, [X_{m+1,
0}^-,[X_{1,0}^-
,\ldots ,[X_{m-1,0}^-,X_{m,1}^-]_{q^{1/2}}\ldots ]_{q^{1/2}}f(k_0),\\
f(x_0^-)&=& \mu f(k_0^{-1})[X_{n,0}^+,[X_{n-1, 0}^+,\ldots ,[X_{m+1,
0}^+,[X_{1,0}^+
,\ldots ,[X_{m-1,0}^+,X_{m,-1}^+]_{q^{1/2}}\ldots ]_{q^{1/2}},\end{eqnarray*}
where $\mu\in\cit$\/ is determined by
$$[f(x_0^+), f(x_0^-)] = \frac{f(k_0) - f(k_0^{-1})}{q-q^{-1}}.\;\;\;\;\Box$$

\eep

\brk Using the relations in ${\cal A}$,  it is not difficult to see that the
isomorphism $f$\/ is independent of the choice of $m\in\{1,2,\ldots , n\}$.
\erk

\subsection{} The following result is proved in [2], [3].
\bep Let $W$\/ be a finite--dimensional irreducible \uqslp--module of type 1.
Then,

\noindent (a) $W$\/ is generated by a vector $w_0$\/ satisfying

$$X_{i,r}^+. w_0 = 0, \;\; \Phi_{i,r}^{{}\pm{}}.w_0 =\phi_{i,r}^{{}\pm{}}w_0 $$
for all $i\in\{1,\ldots ,n\}$, $r\in{\bf Z}$, and some
$\phi_{i,r}^{{}\pm{}}\in\ci$.

\noindent (b) There exist unique monic polynomials $P_1(u),\ldots, P_n(u)$\/
(depending on $W$)  such that the $\phi_{i,r}^{{}\pm{}}$\/ satisfy
$$\sum_{r=0}^{\infty}\phi_{i,r}^+u^r =q^{{\rm
deg}\:P_i}\frac{P_i(q^{-2}u)}{P_i(u)} =\sum_{r=0}^{\infty}\phi_{i,r}^-u^{-r},$$
in the sense that the left and right-hand sides  are the Laurent expansions of
the middle term about $0$\/ and $\infty$\/ respectively.
Assigning to $W$\/ the corresponding $n$--tuple of polynomials defines a one to
one correspondence between the isomorphism classes of finite--dimensional
irreducible \uqslp--modules of type 1 and the set of $n$--tuples of monic
polynomials in one variable $u$. $\Box$
\eep

A consequence of this proposition is:
\becor Let $W$\/ be a finite-dimensional irreducible representation of \uqsla\/
with associated polynomials $P_i$. Set $\lambda =({\rm deg}\;P_1,\ldots, {\rm
deg}\; P_n)$. Then $W$\/ contains the irreducible \uqsl--module  $V(\lambda)$\/
with multiplicity one. Further, if $V(\mu)$\/ is any other \uqsl--module
occurring in $W$, then $\lambda\ge\mu$. $\Box$
\eecor

\subsection{} The next proposition can be proved by studying the action of the
comultiplication $\Delta$\/ of \uqsla\/  on the generators $X_{i,r}^+$\/ etc.,
as in [2].
\bep  Let $W$\/ and $W'$\/ be two finite--dimensional irreducible
\uqsla--modules with associated monic polynomials $P_i$\/ and $P_i'$,
$i=1,\ldots,n$. Let $w_0$\/ and $w_0'$\/ be the generating vectors of $W$\/ and
$W'$\/ as in Proposition 6.2. Then, in $W\ot W'$\/ we have
$$X_{i,r}^+. (w_0\ot w_0') = 0$$ for all $i\in\{1,\ldots ,n\}$, $r\in{\bf Z}$.
Further, $w_0\ot w_0'$\/ is a common eigenvector of the
$\Phi_{i,r}^{{}\pm{}}$\/ with eigenvalues given as in Proposition 6.2 (b) by
the polynomials $P_iP_i'$. $\Box$

\eep This result suggests the following
\bed If $i\in\{1,\ldots ,n\}$, $a\in\cit$, the irreducible finite--dimensional
representation of \uqsla\/ with associated polynomials
$$ P_j(u) =\left\{\begin{array}{ll} u-a\; \;& {\rm if}\; j =i,\\
1\;\; & {\rm otherwise,} \end{array}\right.$$
is called the $i^{th}$\/ fundamental representation of \uqsla\/ with parameter
$a$, and is denoted by $V(\lambda_i, a)$.
\eed
\brk Note that it follows from Corollary 6.2 that $V(\lambda_i ,a)\cong
V(\lambda_i)$\/ as \uqsl--modules.
\erk

\subsection{} We shall need the following result in Section 7.
\bel Let $v_{\lambda_m}$\/ be the  \uqsl--highest weight vector in $V(\lambda_m
,a)$,  where $m\in\{1,\ldots ,n\}$, $a\in\cit$. Then,
$$x_0^+.v_{\lambda_m} = (-1)^{m-1}a^{-1}x_n^-x_{n-1}^-\ldots
x_{m+1}^-x_1^-\ldots x_m^-. v_{\lambda_m}.$$
\eel
Proof. By Proposition 2.3 and the preceding remark,  we know that the weight
spaces of $V(\lambda_m ,a)$\/ as a \uqsl--module are all one--dimensional and
that the  weights are precisely $\epsilon_{i_1}+\epsilon_{i_2}+\ldots
+\epsilon_{i_m}$, \/ $1\le i_1<i_2<\ldots <i_m\le n+1$.
It follows that
$$X^-_{m,1} .v_{\lambda_m} = b x_m^- .v_{\lambda_m}$$
for some $b\in \ci$. Using Proposition 6.1 we get
$$\Phi_{m,1}^+.v_{\lambda_m} = b(q-q^{-1})v_{\lambda_m}.$$
Hence, from Proposition 6.2 (b), we get
$$q(q^{-2}u-a) = (u-a)(q+b(q-q^{-1})u +O(u^2)),$$
so that $b =a^{-1}$. Finally, from Proposition 6.1 again, we find that
$$x_0^+.v_{\lambda_m} = (-1)^{m-1} a^{-1}x_n^-x_{n-1}^-\ldots
x_{m+1}^-x_1^-\ldots x_m^-. v_{\lambda_m}.\;\;\Box$$

\section{Comparison with results of Zelevinsky and Rogawski}
In this  section, we describe a parametrization, due to Zelevinsky [13] and
Rogawski [12], of the finite--dimensional irreducible \hma--modules. We then
relate this, via the functor ${\cal F}$\/ defined in Theorem 4.2, to the
parametrization of the finite--dimensional irreducible \uqsla--modules given in
Section 6.
\subsection{} Since $q$\/ is not  a root of unity, $\hm\cong \ci[\sg]$\/ as an
algebra. It follows that the finite--dimensional \hm--modules are completely
reducible and that the irreducibles are in one to one correspondence with the
partitions of $\ell$.  We now describe this correspondence.

The defining relations of
\hm\/ imply that, if $w\in\sg$\/ and if
$$w =\tau_{i_1}\tau_{i_2}\ldots\tau_{i_k}$$
is any reduced expression for $w$\/ in terms of the simple transpositions
$\tau_i = (i,i+1)$, the element
$$\sigma_w =\sigma_{i_1}\sigma_{i_2}\ldots\sigma_{i_k}\in\hm$$
depends only on $w$.

Let $\le$\/ be the Bruhat order on \sg , and for $w'\le w$, let $P_{w',w}(q)$\/
be the Kazhdan--Lusztig polynomial (see [9]). Define elements $C_w\in\hm$\/ by
$$C_w =q^{\ell(w)}\sum_{w'\le w}(-1)^{\ell(w)
-\ell(w')}q^{-2\ell(w')}P_{w,w'}(q^{-2})\sigma_w .$$
We write $C_i$\/ for $C_{\tau_i}$. Note that $C_i  =q^{-1}\sigma_i -q$. It is
known (see [9]) that $\{C_w\}_{w\in W}$\/ is a basis of \hm, and that
\begin{equation}
C_w\sigma_i = -C_w\; {\rm if}\; w\tau_i < w.
\end{equation}

Let $\ell =\ell_1+\ell_2+\cdots +\ell_p$\/ be a partition $\pi$\/ of $\ell$,
with each $\ell_r >0$, and let
$S_{\ell}^{\pi}$\/ be the subgroup $S_{\ell_1}\times
S_{\ell_2}\times\cdots\times S_{\ell_p}$\/ of \sg\/ which fixes $\pi$.  Let
$w_{r}$\/ be the longest element of the subgroup $S_{\ell_r}$, i.e. the
permutation which reverses the order of
$(\ell_1+\ell_2+\cdots +\ell_{r-1}+1,\ldots ,\ell_1+\cdots +\ell_r)$,
 and set
$w_{\pi} =w_1w_2\ldots w_p$.
Let $I_{\pi}$\/ be the right ideal in \hm\/ generated by $C_{w_{\pi}}$.
\bep {{\rm ([12])}} For every partition $\pi$\/ of $\ell$, $I_{\pi}$\/ has a
unique irreducible quotient $J_{\pi}$\/ in which $C_{w_{\pi}}$\/ has non--zero
image. Conversely, every finite--dimensional irreducible right \hm-module is
isomorphic to some $J_{\pi}$. $\Box$\eep

\subsection{} Using Jimbo's functor ${\cal J}$, we can compare this
parametrization of the finite--dimensional irreducible representations of \hm\/
with that of the representations of \uqsl\/ given by their highest weights.
\bep Let $1\le\ell \le n$\/ and let $\ell_1+\ell_2+\cdots +\ell_p$\/ be a
partition $\pi$\/ of $\ell$. Then,
$${\cal J}(J_\pi)\cong V(\lambda_{\ell_1}+\lambda_{\ell_2}+\cdots
+\lambda_{\ell_p})$$
as \uqsl--modules.
\eep
Proof. We need the following lemma, which follows from (1).
\bel Let $\pi$\/ be as in the prceding proposition, and let $1\le i\le \ell$\/
be such that $i\ne\sum_{j=1}^r\ell_j$\/ for any $1\le r<p$. Let $\unv\in\vm$\/
have $v_r\ot v_s$\/ in the $i^{th}$\/ and $(i+1)^{th}$\/ positions, and let
${\bf v}'$\/ be the result of interchanging the vectors in these positions.
Then, in ${\cal J}(J_{\pi})$, we have
$$C_{w_{\pi}}\ot\unv' =\left\{\begin{array}{ll} -q^{-1}C_{w_{\pi}}\ot\unv& {\rm
if}\; r< s,\\
-qC_{w_{\pi}}\ot\unv&{\rm if}\; r>s,\\
0&{\rm if}\; r=s.\end{array}\right.\;\;\;\Box$$
\eel

Returning to the proof of the proposition, note that the weight space of \vm\/
of weight $\lambda_{\ell_1}+\lambda_{\ell_2}+\cdots +\lambda_{\ell_p}$\/ is
spanned by the permutations of the vector
$$\unv_{\pi} =v_1\ot v_2\cdots\ot v_{\ell_1}\ot v_1\ot v_2\ot\cdots\ot
v_{\ell_2}\ot v_1\cdots\ot v_{\ell_p} .$$
By Proposition 4.1, there exists a partition $\pi'$\/ of $\ell$, say $\ell =
\ell_1'+\ell_2'+\cdots + \ell_r'$, such that
\begin{equation}
{\cal J}(J_{\pi'})\cong V(\lambda_{\ell_1}+\cdots +\lambda_{\ell_r}).
\end{equation}
By the lemma, if $v_{i_1}\ot\cdots\ot v_{i_\ell}$\/ is any permutation of ${\bf
v}_{\pi}$,
$$C_{w_{\pi}}\ot v_{i_1}\ot\cdots\ot v_{i_\ell} = 0$$
unless the first $\ell_1'$\/ vectors in the sequence $v_{i_1},\ldots ,
v_{i_{\ell}}$\/ are distinct,  together with the next $\ell_2'$, $\ldots$, and
the last $\ell_r'$. It follows that, if $\le$\/ is the usual lexicographic
ordering on the set
of partitions of $\ell$, we have $\pi'\le\pi$. But the map $\pi\to\pi'$\/
defined by (13) is  a bijection since ${\cal J}$\/ is  an equivalence. Since
$\le$\/ is a total ordering it follows that this bijection is
the identity map, i.e. $\pi'=\pi$. $\Box$
\subsection{}
We now turn to the representations of affine Hecke algebras. Recall the
universal modules $M_{\bf a}$\/ defined in Section 3.4. We begin with the
following elementary result.
\bel{{\rm ([12])}} Let ${\bf a}=(a_1,\ldots,a_\ell)\in(\ci^\times)^\ell$, $w\in
S_\ell$, $j\in\{1,\ldots,\ell\}$. Then, in $M_{\una}$, we have
$$C_w.y_j = a_{w^{-1}(j)} C_w +\sum_{w'<w}\alpha_{w'}C_{w'}$$
for some $\alpha_{w'}\in\ci$. $\Box$
\eel

\subsection{}Following Rogawski [12] and Zelevinsky [13], we make the following
definition.
\bed The segment $s$\/ with centre $a\in\ci^\times$\/ and length $|s|=k$\/ is
the ordered sequence $s=(aq^{-k+1},aq^{-k+3},\ldots,aq^{k-1})\in(\cit)^k$.
\eed

If $\uns =\{s_1, s_2,\ldots ,s_p\}$\/ is any (unordered) collection of
segments, and if $|s_r| =\ell_r$, then
$\ell =\ell_1+\ell_2+\cdots +\ell_p$\/ is a partition $\pi(\uns)$\/ of $\ell$.
\bep{{\rm ([12])}} Let $\ell\ge 1$\/ and let ${\bf s}=\{s_1,\ldots,s_p\}$\/ be
any collection of segments, the sum of whose lengths is $\ell$. Let $\una
=(s_1,\ldots ,s_p)\in(\cit)^{\ell}$\/ be the result of juxtaposing the segments
in \uns . Then,

\noindent (a) $I_{\pi({\bf s})}$\/ is an \hma--submodule of $M_{\bf a}$ (this
statement makes sense in view of Proposition 3.4 (b));

\noindent (b) with the \hma--module structure from $M_{\bf a}$, $I_{\pi({\bf
s})}$\/ has a unique irreducible subquotient $V_{\bf a}$\/ in which
$C_{w_{\pi({\bf s})}}$\/ has non--zero image.

Moreover, every finite--dimensional irreducible right \hma--module  is
isomorphic to some $V_{\bf a}$. $\Box$
\eep
\subsection{} To prove the main result of this section, we shall need another
description of $I_{\pi({\bf s})}$\/ (we continue to use the notation of Section
7.4). Let $\Sigma^{\pi({\bf s})}\subset\sg$\/ be  the set of transpositions
$\tau_i =(i,i+1)$\/ for $i\in\{1,\ldots
,\ell\}\backslash\{\ell_1,\ell_1+\ell_2,\ldots,\ell_1+\cdots+\ell_{p-1}\}$. For
$\tau_i\in\Sigma^{\pi({\bf s})}$, let ${\bf a}_{\tau_i}$\/ be the result of
interchanging
the $i^{th}$\/ and $(i+1)^{th}$\/ components of \una , and let
$$A_{{\bf a} ,i}:M_{{\bf a}_{\tau_i}}\to M_{\bf a}$$
be the map given by left multiplication by $C_i$ (we identify $M_{\bf a}$\/ and
$M_{{\bf a}_{\tau_i}}$\/ with \hma\/ in the usual way).
\bep {\rm ([12])}
With the above notation:

\noindent (a) $A_{{\bf a},i}$\/ is  a homomorphism of \hma--modules;

\noindent (b) regarded as an \hma--submodule of $M_{\bf a}$,
$$I_{\pi({\bf s})} = \bigcap_{\tau_i\in\Sigma^{\pi({\bf s})}}({\rm image\;
of}\; A_{{\bf a} ,i}).\;\;\Box$$
\eep
\subsection{} We can now state the main result of this section.
\bet
Let $\uns = \{s_1,\ldots, s_p\}$\/ be a collection of segments, the sum of
whose lengths is $\ell$, let $a_r$\/ be the centre of $s_r$\/ and $\ell_r$\/
its length, and let ${\bf a}=(s_1,\ldots,s_p)\in(\ci^\times)^\ell$\/ be the
result of juxtaposing $s_1,\ldots,s_p$, as in Proposition 7.4.
Then, if $\ell\le n$, ${\cal F}(V_{\bf a})$\/ is the irreducible \uqsla--module
defined by the polynomials
$$P_i(u) =\prod_{\{j|\ell_j=i\}}(u-a_j^{-1}), \;\;i=1,\ldots ,n.$$
\eet
Proof. We first prove the result in the special case $p =1$, so that $\una
=(aq^{-\ell+1}, aq^{-\ell+3},\ldots ,aq^{\ell-1})$\/ (we drop the subscripts
for simplicity). Note that $w_{\pi({\bf s})} = w_0$, the longest element of \sg
, and that $I_{\pi({\bf s})}$ (= $J_{\pi({\bf s})} = V_{\bf a}$)\/ is
one-dimensional and spanned by
$C_{w_0}$. By Proposition 7.2,
$${\cal J}(I_{\pi({\bf s})})\cong V(\lambda_{\ell}),$$
the highest weight vector being
$$\unv_{\lambda_\ell} = C_{w_0}\ot v_1\ot v_2\ot\cdots \ot v_{\ell} .$$
As a \uqsla--module, ${\cal F}(V_{\bf a})$\/ is therefore defined by the
polynomials
$$P_i(u) =\left\{\begin{array}{ll} u -a' &{\rm if}\; i =\ell,\\
1&{\rm otherwise,}\end{array}\right.$$
for some $a'\in\cit$. To compute $a'$, note first that, by the definition of
${\cal F}$,
$$x_0^+.\unv_{\lambda_\ell} =C_{w_0}.y_1\ot v_{n+1} \ot v_2\ot\cdots\ot
v_{\ell} .$$
Since $I_{\pi({\bf s})}$\/ is one--dimensional, Lemma 7.3 implies that
\begin{equation}x_0^+.\unv_{\lambda_\ell} =q^{\ell-1}a C_{w_0}\ot v_{n+1}\ot
v_2\ot\cdots\ot v_{\ell}.\end{equation}
On the other hand Lemma 6.4 gives
\begin{eqnarray*}x_0^+.\unv_{\lambda_\ell} &=&
(-1)^{\ell-1}(a')^{-1}x_n^-x_{n-1}^-\cdots x_{\ell+1}^-x_1^-x_2^-\cdots
x_{\ell}^-.\unv_{\lambda_\ell}\\
&=& (-1)^{\ell-1} (a')^{-1}(C_{w_0}\ot v_2\ot\cdots\ot v_{\ell}\ot
v_{n+1}).\end{eqnarray*}
Now by (12),
$$C_{w_0}\sigma_i^{-1} = -C_{w_0},$$
and by (1),
$$v_r\ot v_{n+1} = q\check{R}^{-1}(v_{n+2}\ot v_r),\;\;{\rm if}\; r\le n.$$
Hence,
$$C_{w_0}\ot v_2\ot\cdots\ot v_{\ell}\ot v_{n+1} =
(-1)^{\ell-1}q^{\ell-1}C_{w_0}\ot v_{n+1}\ot v_2\ot\cdots\ot v_{\ell},$$
and so
$$x_0^+.\unv_{\lambda_\ell} = q^{\ell-1}(a')^{-1}C_{w_0}\ot v_{n+1}\ot
v_2\ot\cdots\ot v_{\ell}.$$
Comparing with (14) gives $a'=a^{-1}$. (It follows from the proof of
Proposition 7.2 that $C_{w_0}\ot v_{n+1}\ot v_2\ot\cdots\ot v_{\ell}\ne 0$.)

Suppose now that $r$\/ is arbitrary. From Proposition 7.5 (b),
\begin{equation}
{\cal F}(I_{\pi({\bf s})}) =\bigcap_{\tau_i\in\Sigma^{\pi({\bf s})}}( {\rm
image\;of\;} {\cal F}(A_{{\bf a},i})).
\end{equation}
To compute ${\cal F}(A_{{\bf a},i})$, note that $\unv\mapsto 1\ot\unv$\/
defines an isomorphism of \uqsl--modules $V^{\ot\ell}\to{\cal F}(M_{\bf a})$,
and that
$${\cal F}(A_{{\bf a},i})(1\ot\unv)=C_i\ot\unv=1\ot C_i.\unv.$$
It follows that
$${\cal F}(A_{{\bf a} ,i}) = q^{-1}\check{R}_i - q\in {\rm End}_{\ci}\;(\vm
).$$
{}From (15) and the $r=1$\/ case, it follows that
$${\cal F}(I_{\pi({\bf s})}) = V(\lambda_{{\ell_1}},a_1^{-1})\ot\cdots\ot
V(\lambda_{\ell_p},a_p^{-1}). $$
By Propositions 7.2 and 7.4 (b), ${\cal F}(V_{\bf a})$\/ is the unique
irreducible subquotient of ${\cal F}(I_{\pi({\bf s})})$\/ in which the tensor
product of the highest weight vectors in the $V(\lambda_{\ell_r}, a_r^{-1})$\/
has non--zero image. The theorem now follows from the multiplicativity of the
polynomials in Proposition 6.3. $\Box$

\pagebreak

\vskip 36pt
\noindent {\bf{ \Large References}}
\vskip 24pt
\noindent 1. Chari, V. and Pressley, A.N., New unitary representations of loop
groups, Math. Ann. 275 (1986), 87--104.

\noindent 2. Chari V. and Pressley A.N., Quantum affine algebras, Comm. Math.
Phys. 142 (1991), 261--283.

\noindent 3. Chari V. and Pressley A.N., {\it A Guide to Quantum Groups},
Cambridge University Press, Cambridge (in press).

\noindent 4. Cherednik, I. V., A new interpretation of Gelfand-Tzetlin bases,
Duke Math. J. 54 (1987), 563--577.

\noindent 5. Drinfeld, V. G., Degenerate affine Hecke algebras and Yangians,
Func. Anal. Appl. 20 (1986), 62--64.

\noindent 6. Drinfeld, V. G., A new realization of Yangians and quantized
affine algebras, Sov. Math. Dokl. 36 (1988), 212--216.

\noindent 7. Jimbo, M., A $q$--analogue of $U(gl(N+1))$, Hecke algebra and the
Yang--Baxter equation, Lett. Math. Phys. 11 (1986), 247--252.

\noindent 8. Kac, V. G., {\it Infinite Dimensional Lie Algebras}, Birkh\"auser,
Boston, 1983.

\noindent 9. Kazhdan, D. and Lusztig, G., Representations of Coxeter groups and
Hecke algebras, Invent. Math. 53 (1979), 165--184.

\noindent 10. Lusztig, G., Quantum deformations of certain simple modules over
enveloping algebras, Adv. Math. 70 (1988), 237-249.

\noindent 11. MacLane, S., {\it Categories for the Working Mathematician},
Springer, 1971.

\noindent 12. Rogawski, J.D., On modules over the Hecke algebra of a $p$--adic
group, Invent. Math. 79 (1985), 443--465.

\noindent 13. Zelevinsky, A.V., Induced representations of reductive $p$--adic
groups II. On irreducible representations of $GL_n$, Ann. Sci. Ec. Norm. Sup.
$4^e$\/ S\'er. 13 (1980), 165--210.
\vskip 24pt

\noindent Vyjayanthi Chari,

\noindent Department of Mathematics,

\noindent University of California,

\noindent Riverside, CA 92521, U.S.A.

\vskip 12pt

\noindent Andrew Pressley,

\noindent Department of Mathematics,

\noindent King's College, Strand,

\noindent London WC2R 2LS, U.K.
\end{document}